\documentclass[12pt]{article}
\usepackage{amsfonts}
\usepackage{tabularx}
\begin{document}

\thispagestyle{plain}

\begin{center}
\Large{Mahler measure of some $n$-variable polynomial families}

\medskip
\normalsize{Matilde N. Lal\'{\i}n
{\footnote[1]{E-mail address: mlalin@math.utexas.edu}}
}

\medskip
University of Texas at Austin. Department of Mathematics.
 1 University Station C1200.  Austin, TX 78712, USA
\end{center}



\begin{abstract}
The Mahler measures of some $n$-variable polynomial families are given in terms
of special values of the Riemann zeta function and a Dirichlet L-series,
generalizing the results of  \cite{L}. The technique introduced in this work
also motivates certain identities among Bernoulli numbers and symmetric
functions. \end{abstract}

\newcommand{\Cset}{\mathbb{C}}
\newcommand{\F}{\mathbb{F}}
\newcommand{\Rset}{\mathbb{R}}
\newcommand{\Qset}{\mathbb{Q}}
\newcommand{\Nset}{\mathbb{N}}
\newcommand{\Zset}{\mathbb{Z}}
\newcommand{\PP}{\mathbb{P}}
\newcommand{\HH}{\mathbb{H}}
\newcommand{\MU}{\mathbb{\mu}}
\newcommand{\TT}{\mathbb{T}}
\newcommand{\binom}[2]{{#1\choose#2}} 
\newcommand{\Res}{\mathrm{Res}}
\newcommand{\Li}{\mathrm{Li}} 
\newcommand{\I}{\mathrm{I}} 
\newcommand{\re}{\mathop{\mathrm{Re}}} 
\newcommand{\im}{\mathop{\mathrm{Im}}} 
\newcommand{\ii}{\mathrm{i}} 
\newcommand{\Lf}{\mathrm{L}} 
\newcommand{\dd}{\mathrm{d}}
\newcommand{\e}{\mathrm{e}}
\newcommand{\mod}{\mathrm{mod}}
\newcommand{\pf}{\noindent {\bf PROOF.} \quad}
\newcommand{\qed}{$\Box$}
\newtheorem{thm}{Theorem}
\newtheorem{defn}[thm]{Definition}
\newtheorem{prop}[thm]{Proposition}
\newtheorem{cor}[thm]{Corollary}
\newtheorem{lem}[thm]{Lemma}
\newtheorem{conj}[thm]{Conjecture}
\newtheorem{rem}[thm]{Remark}
\newtheorem{ex}[thm]{Example}
\newtheorem{exs}[thm]{Examples}
\newtheorem{obs}[thm]{Observation}
\section{Introduction}

The goal of this work is to exhibit three families of multivariable polynomials whose Mahler measure depend (in most of the cases) on special values of the Riemann zeta function and the L-series on the Dirichlet character of conductor four.

The Mahler measure of a polynomial $P \in \Cset[x_1, \dots x_n]$ is defined as
\[m(P) = \frac{1}{(2\pi \ii)^n} \int_{\TT^n} \log | P(x_1, \dots, x_n)| \frac{\dd 
x_1}{x_1} \dots \frac{\dd x_n}{x_n} .\] Here $\TT^n = \{ (z_1, \dots, z_n) \in
\Cset ^n \,|\, |z_1| = \dots = |z_n| = 1 \}$ is the unit torus.

For one-variable polynomials, Jensen's formula gives a simple expression for the Mahler measure as a function on the roots of the polynomial. However,  it is in general  a very hard problem to give an explicit closed formula for the Mahler measure of a polynomial in two or more variables.

For up to five variables, several examples have been produced by Bertin
\cite{Be}, Boyd  \cite{B1,B2,B3}, Boyd and Rodriguez-Villegas \cite{BRV1,BRV2},
Cassaigne and Maillot \cite{M}, Condon \cite{C}, Smyth \cite{S1,S2}, Vandervelde
\cite{V}, the author \cite{L}, among others.

Smyth \cite{S3} gave and example of an $n$-variable family of polynomials whose
Mahler measures are related to special values of hypergeometric series.

We have analyzed the $n$-variable versions of the polynomials studied in
\cite{L} and found closed formulas for their Mahler measures, which in most of
the cases depend on special values of the Riemann zeta function and Dirichlet
L-series. More precisely, we have proved that

\begin{thm} We have the following identities: \footnote{In order to
simplify notation, we describe the polynomials as rational functions, writing
$1+\frac{1-x}{1+x}z$ instead of $1+x+(1-x)z$, and so on. The Mahler measure does
not change since the denominators are product of cyclotomic polynomials.}

(i) For $n \geq 1$:
 \[ \pi^{2n} m\left( 1+  \left( \frac{1-x_1}{1+x_1}\right) \dots
\left( \frac{1-x_{2n}}{1+x_{2n}}\right) z \right) \]
\[=\frac{1}{(2n-1)!}\sum_{h=1}^n  s_{n-h}(2^2, \dots, (2n-2)^2)
\frac{(2h)!(2^{2h+1}-1)}{2}\pi^{2n-2h} \zeta(2h+1).\]
 \begin{equation}
\label{eq1}
\end{equation}

For $n \geq 0$:
\[ \pi^{2n+1} m\left( 1+  \left( \frac{1-x_1}{1+x_1}\right) \dots
\left( \frac{1-x_{2n+1}}{1+x_{2n+1}}\right) z \right)\]
\[= \frac{1}{(2n)!} \sum_{h=0}^n s_{n-h}(1^2, \dots, (2n-1)^2)  (2h+1)! 2^{2h+1}
\pi^{2n-2h} \Lf(\chi_{-4}, 2h+2).\]
\begin{equation}\label{eq2}
\end{equation}

(ii) For $n\geq 1$:
\[ \pi^{2n+2} m\left(  1+x +  \left( \frac{1-x_1}{1+x_1}\right)
\dots \left( \frac{1-x_{2n}}{1+x_{2n}}\right)(1+y) z
\right)\]
\[=\frac{1}{(2n-1)!}\sum_{h=1}^n \frac{(2h+2)! (2^{2h+3}-1)}{8}\]
\[\left( \sum_{l=0}^{n-h}s_{n-h-l}(2^2, \dots, (2n-2)^2) \binom{2(l+h)}{2h}
(-1)^l \frac{2^{2l}}{l+h} B_{2l}\right) \pi^{2n-2h} \zeta(2h+3). \]
\begin{equation}\label{eq3}
\end{equation}

For $n \geq 0$:
\[ \pi^{2n+3} m\left(  1+x +  \left( \frac{1-x_1}{1+x_1}\right)
\dots \left( \frac{1-x_{2n+1}}{1+x_{2n+1}}\right)(1+y) z \right)\]
\[=\frac{1}{(2n)!} \sum_{h=0}^{n} s_{n-h}(1^2, \dots, (2n-1)^2)  \]
\[2^{2h+1} \pi ^{2n-2h} \left(\ii (2h)!
\mathcal{L}_{3,2h+1}(\ii,\ii) + (2h+1)!\pi^2 \Lf(\chi_{-4}, 2h+2)\right).\]
\begin{equation}\label{eq4}
\end{equation}

(iii) For $n\geq 1$:
\[ \pi^{2n+1} m \left ( 1 + \left( \frac{1 - x_1}{1+x_1}
\right)\dots  \left( \frac{1 - x_{2n}}{1+x_{2n}} \right)x + \left( 1 -\left(
\frac{1 - x_1}{1+x_1} \right)\dots  \left( \frac{1 - x_{2n}}{1+x_{2n}}
\right)\right) y  \right) \] \[= \frac{\pi^{2n+1}}{2} \log 2\]
\[+\frac{1}{(2n-1)!}\sum_{h=1}^n s_{n-h}(2^2, \dots, (2n-2)^2) \frac{
(2h)!(2^{2h+1} - 1)}{4}\pi ^{2n-2h+1}  \zeta(2h+1)\]
 \[+ \frac{1}{(2n-1)!}\sum_{h=1}^n \frac{(2h)!(2^{2h+1}-1)}{4} \]
\[\left(\sum_{l=0}^{n-h}  s_{n-h-l}(2^2, \dots,
(2n-2)^2) \binom{2(l+h)}{2l} (-1)^{l+1}\frac{2^{2l}(2^{2l-1}-1)}{l+h}  B_{2l}
\right)\pi^{2n-2h+1} \zeta(2h+1).\]
\begin{equation}\label{eq5}
\end{equation}

For $n\geq 0$:
\[ \pi^{2n+2} m \left ( 1 + \left( \frac{1 - x_1}{1+x_1}
\right)\dots  \left( \frac{1 - x_{2n+1}}{1+x_{2n+1}} \right)x + \left( 1 -\left(
\frac{1 - x_1}{1+x_1} \right)\dots  \left( \frac{1 - x_{2n+1}}{1+x_{2n+1}}
\right)\right) y  \right) \] \[= \frac{\pi^{2n+2}}{2} \log 2\]\[ +
\frac{1}{(2n+1)!}\sum_{h=0}^{n} s_{n-h}(2^2,\dots,(2n)^2)      \frac{(2h+2)!
(2^{2h+3} - 1)}{4}  \pi^{2n-2h} \zeta(2h+3)\] \[+ \frac{1}{(2n-1)!}\sum_{h=1
}^{n} \frac{(2h)! (2^{2h+1}-1)}{4 } \]
\[\left(\sum_{l=0}^{n-h}  s_{n-h-l}(2^2, \dots, (2n-2)^2) \binom{2(l+h)}{2l}
(-1)^{l+1}\frac{2^{2l}(2^{2l-1}-1)}{l+h}  B_{2l} \right)\pi^{2n-2h+2}
\zeta(2h+1).\]
\begin{equation}\label{eq6}
\end{equation}
\end{thm}

Where $B_h$ is the $h$-Bernoulli number, $\frac{x}{\e^x-1} = \sum_{n=0}^\infty
\frac{B_n x^n}{n!}$.

$\zeta$ is the Riemann zeta function,
\[\Lf(\chi_{-4}, s) := \sum_{n=1}^\infty \frac{\chi_{-4}(n)}{n^s}  \]
\[\chi_{-4}(n) =  \left \{ \begin{array}{cl}
\left(\frac{-1}{n}\right) & \quad \mbox{if}\; n \; \mbox{odd } \\\\
0 & \quad\mbox{if}\; n \; \mbox{even}
\end{array} \right. \]
and $\mathcal{L}_{r,s}(\alpha, \alpha)$ are linear combinations of
multiple polylogarithms (they will be defined later).

Also,
\begin{equation}
s_l(a_1,\dots,a_k) = \left \{ \begin{array}{ll}1 & \mbox{if} \quad l=0 \\
 \sum_{i_1 < \dots < i_l} a_{i_1} \dots a_{i_l} & \mbox{if} \quad 0< l \leq k \\
0 & \mbox{if} \quad k < l
\end{array} \right .
\end{equation}
are the elementary symmetric polynomials, i.e.,
\begin{equation} \label{eq:sym}
\prod_{i=1}^k (x + a_i) = \sum_{l=0}^k s_l(a_1,\dots,a_k) x^{k-l}.
\end{equation}

For concreteness, we list the first values for each family in the following tables:

\begin{center}
\renewcommand{\arraystretch}{2.7}
\newcolumntype{Y}{>{\centering\arraybackslash}X}%
 \begin{tabularx}{160mm} %
{|>{\setlength{\hsize}{1 \hsize}}Y| >{\setlength{\hsize}{1 \hsize}}Y|}
\hline
$\pi^2 m \left( 1 + \left( \frac{1-x_1}{1+x_1}\right ) \left( \frac{1-x_2}{1+x_2} \right ) z \right) $&  $  7\, \zeta(3)$\\
$\pi^4  m \left( 1 + \left( \frac{1-x_1}{1+x_1}\right )  \dots \left( \frac{1-x_4}{1+x_4}\right )  z \right) $  &  $ 62 \zeta(5) + \frac{14 \pi^2}{3}  \zeta(3)$\\
$\pi^6 m \left( 1 + \left( \frac{1-x_1}{1+x_1}\right )  \dots  \left( \frac{1-x_6}{1+x_6}\right )  z \right) $ & $ 381 \zeta(7) +62 \pi^2  \zeta(5)+ \frac{56 \pi^4}{15}  \zeta(3) $ \\
$\pi^8 m \left( 1 + \left( \frac{1-x_1}{1+x_1}\right )  \dots  \left( \frac{1-x_8}{1+x_8}\right )  z \right) $ & $ 2044 \zeta(9) +508 \pi^2  \zeta(7)+ \frac{868 \pi^4}{15}  \zeta(5) + \frac{16\pi^6}{5} \zeta(3)$ \\

\hline
$\pi m \left( 1 + \left( \frac{1-x_1}{1+x_1}\right )  z \right) $ & $2  \Lf (\chi_{-4}, 2)  $ \\
$\pi^3 m\left(  1 +\left(  \frac{1-x_1}{1+x_1}\right )  \dots  \left( \frac{1-x_3}{1+x_3}\right )  z \right)$ & $24 \Lf (\chi_{-4}, 4) + \pi^2 \Lf (\chi_{-4}, 2)  $ \\
$\pi^5 m\left(1 + \left( \frac{1-x_1}{1+x_1}\right )  \dots \left(  \frac{1-x_5}{1+x_5}\right )   z\right) $&$
160  \Lf(\chi_{-4}, 6) + 20 \pi^2 \Lf(\chi_{-4},4) + \frac{3 \pi^4}{4} \Lf(\chi_{-4},2)$\\
$\pi^7 m\left(1 + \left( \frac{1-x_1}{1+x_1}\right )  \dots \left(  \frac{1-x_7}{1+x_7}\right )   z\right) $&$
896  \Lf(\chi_{-4}, 8) + \frac{560}{3} \pi^2 \Lf(\chi_{-4},6) + \frac{259}{15} \pi^4\Lf(\chi_{-4},4)+\frac{5}{8}\pi^6\Lf(\chi_{-4},2) $\\
\hline
\end{tabularx}
\end{center}

\bigskip

\begin{center}
\renewcommand{\arraystretch}{2.7}
\newcolumntype{Y}{>{\centering\arraybackslash}X}%
 \begin{tabularx}{160mm} %
{|>{\setlength{\hsize}{1 \hsize}}Y| >{\setlength{\hsize}{1 \hsize}}Y|}
\hline
$ \pi^2 m\left( 1+x  + (1+y)z \right) $& $ \frac{7}{2} \zeta(3)$\\
$ \pi^4 m\left( 1+x  +\left(  \frac{1-x_1}{1+x_1} \right ) \left( \frac{1-x_2}{1+x_2}\right ) (1+y)z \right) $& $ 93 \zeta(5)$\\
$ \pi^6 m\left( 1+x  +\left(  \frac{1-x_1}{1+x_1}\right ) \dots \left( \frac{1-x_4}{1+x_4}\right ) (1+y)z \right) $& $ \frac{1905}{2} \zeta(7) + 31 \pi^2 \zeta(5)$ \\
$ \pi^8m\left(  1+x  +\left(  \frac{1-x_1}{1+x_1}\right ) \dots \left( \frac{1-x_6}{1+x_6}\right ) (1+y)z \right)   $&$ 7154 \zeta(9)+ 635 \pi^2  \zeta(7) + \frac{248 \pi^4 }{15} \zeta(5)  $\\
\hline
$ \pi^3 m\left( 1+x  +\left(  \frac{1-x_1}{1+x_1}\right ) (1+y)z \right) $&$2 \pi^2 \Lf (\chi_{-4}, 2) +2 \ii  \mathcal{L}_{3,1}(\ii,\ii)$\\
$ \pi^5 m\left( 1+x  +\left(  \frac{1-x_1}{1+x_1}\right ) \dots \left( \frac{1-x_3}{1+x_3}\right ) (1+y)z \right) $&$24 \pi^2 \Lf(\chi_{-4}, 4) + \pi^4 \Lf (\chi_{-4}, 2) +16 \ii \mathcal{L}_{3,3}(\ii,\ii) + 4 \pi \ii \mathcal{L}_{3,1}(\ii,\ii)$\\
\hline
\end{tabularx}
\end{center}

\bigskip

\begin{center}
\renewcommand{\arraystretch}{2.7}
\newcolumntype{Y}{>{\centering\arraybackslash}X}%
 \begin{tabularx}{160mm} %
{|>{\setlength{\hsize}{1.4 \hsize}}Y| >{\setlength{\hsize}{.6 \hsize}}Y|}
\hline
$\pi^3 m \left( 1 + \left( \frac{1-x_1}{1+x_1}\right)\left( \frac{1-x_2}{1+x_2}\right ) x +  \left( 1 - \left( \frac{1-x_1}{1+x_1}\right)\left( \frac{1-x_2}{1+x_2}\right ) \right) y \right) $&  $  \frac{21 \pi}{4}  \zeta(3) + \frac{\pi^3}{2} \log 2 $\\
$\pi^5 m \left( 1 + \left( \frac{1-x_1}{1+x_1}\right)\dots \left( \frac{1-x_4}{1+x_4}\right ) x +  \left( 1 - \left( \frac{1-x_1}{1+x_1}\right)\dots \left( \frac{1-x_4}{1+x_4}\right ) \right) y \right) $&  $  \frac{155 \pi}{4}  \zeta(5) + \frac{14\pi^3}{3} \zeta(3) +\frac{\pi^5}{2} \log 2 $\\
\hline
$\pi^2 m \left( 1 + \left( \frac{1-x_1}{1+x_1}\right ) x +  \left( 1 - \left( \frac{1-x_1}{1+x_1}\right) \right) y \right) $&  $  \frac{7 }{2}  \zeta(3) + \frac{\pi^2}{2} \log 2 $\\
$\pi^4 m \left( 1 + \left( \frac{1-x_1}{1+x_1}\right) \dots  \left( \frac{1-x_3}{1+x_3}\right ) x +  \left( 1 - \left( \frac{1-x_1}{1+x_1}\right)\dots  \left( \frac{1-x_3}{1+x_3}\right ) \right) y \right) $&  $ 31 \zeta(5) + \frac{7 \pi^2}{3}  \zeta(3) + \frac{\pi^4}{2} \log 2 $\\
\hline
\end{tabularx}
\end{center}

\bigskip

\section{An important integral}

Before proving our main result, we will need to prove some auxiliary statements.

We will need to compute the integral $\int_0^\infty \frac{x \log^k x \dd x}{(x^2+a^2)(x^2+b^2)} $. The following Lemma will help:
\begin{lem}\label{lem1} We have the following integral:
\begin{equation} \int_0^\infty \frac{ x^\alpha \dd x}{(x^2+a^2)(x^2+b^2)}=
\frac{\pi (a^{\alpha-1}-b^{\alpha-1})}{2 \cos{\frac{\pi \alpha}{2}} (b^2-a^2)}
\qquad \mathrm{for} \quad 0<\alpha < 1. \end{equation}
\end{lem}

\pf We write the integral as a difference of two integrals:
\begin{equation} \int_0^\infty \frac{x^\alpha \dd x}{(x^2+a^2)(x^2+b^2)}=
\int_0^\infty \left( \frac{1}{x^2+a^2} - \frac{1}{x^2+b^2}\right)
\frac{x^\alpha \dd x}{b^2-a^2}.
\end{equation}
Now, when $0<\alpha<1$,
\[\int_0^\infty \frac{x^\alpha \dd x}{x^2+a^2} = \frac{1}{1-\e^{2 \pi \ii
\alpha}} 2 \pi \ii \sum_{x \not = 0} \Res \left \{ \frac{x^\alpha}{x^2+a^2}\right \} \]
(see, for instance, section 5.3 in chapter 4 of the Complex Analysis book by Ahlfors \cite{A} ). Then,
\[\int_0^\infty \frac{x^\alpha \dd x}{x^2+a^2} =\frac{\pi a^{\alpha -1}}{2 \cos
\frac{\pi \alpha}{2}}.\]
Thus, we get the result.\qed

By continuity, the formula in the statement is true for $\alpha=1$, in fact the integral converges for $\alpha<3$.

Next, we will define some polynomials  that will be used in the formula for $\int_0^\infty \frac{x \log^k x \dd x}{(x^2+a^2)(x^2+b^2)} $.
\begin{defn} Let $P_k(x) \in \Qset[x]$, $k\geq 0$,  be defined recursively as follows:
\begin{equation} \label{defiP}
P_k(x)= \frac{x^{k+1}}{k+1} + \frac{1}{k+1} \sum_{j>1\,\mathrm{(odd )}
}^{k+1}(-1)^{\frac{j+1}{2}} \binom{k+1}{j} P_{k+1-j}(x). \end{equation}
\end{defn}
For instance, the first $P_k(x)$ are:
\begin{eqnarray*}
P_0(x) & = & x\\
P_1(x) & = & \frac{x^2}{2} \\
P_2(x) & = & \frac{x^3}{3} +\frac{x}{3} \\
P_3(x) & = & \frac{x^4}{4} + \frac{x^2}{2}\\
P_4(x) & = &  \frac{x^5}{5} + \frac{2 x^3}{3}  + \frac{7x}{15}\\
P_5(x) & = &  \frac{x^6}{6} + \frac{5 x^4}{6} + \frac{7x^2}{6}
\end{eqnarray*}

\begin{lem}\label{easyprop} The following properties are true
\begin{enumerate}
\item $\deg P_k = k+1$.
\item Every monomial of $P_k(x)$ has degree odd  (even) for $k$ even (odd ).
\item $P_k(0)=0$.
\item $P_{2l}(\ii)=0$ for $l>0$.
\item $(2l+1) P_{2l}(x) = \frac{\partial}{\partial x} P_{2l+1}(x)$.
\item $2l P_{2l-1}(x) \equiv \frac{\partial}{\partial x} P_{2l}(x) \, \mathrm{mod}\, x$.
\end{enumerate}
\end{lem}
The above properties can be easily proved by induction. These properties, together with $P_0$, determine the whole family of polynomials $P_k$ because of the recursive nature of the definition. At this point, it should be noted that this family  is closely related to Bernoulli polynomials. We postpone the discussion of this topic for the Appendix, since the explicit form of the polynomials $P_k$ is barely needed in order to perform the computation of the Mahler measures.

We are now ready to prove the key Proposition for the main Theorem:
\begin{prop}\label{prop} We have:
\begin{equation}
\int_0^\infty \frac{x \log^k x \dd x}{(x^2+a^2)(x^2+b^2)} =
\left(\frac{\pi}{2}\right)^{k+1}\frac{ P_{k}\left(\frac{2\log a}{\pi}\right) -
P_{k}\left(\frac{2\log b}{\pi}\right)}{a^2-b^2}. \end{equation} \end{prop}

\pf The idea, suggested by Rodriguez-Villegas, is to obtain the value of the integral in the statement by differentiating $k$ times the integral of Lemma \ref{lem1} and then evaluating at $\alpha=1$. Let
\[ f(\alpha) = \frac{\pi (a^{\alpha-1}-b^{\alpha-1})}{2 \cos{\frac{\pi \alpha}{2}} (b^2-a^2)}\]
which is the value of the integral in the Lemma \ref{lem1}.  In other words, we have
\[f^{(k)}(1) =\int_0^\infty \frac{x \log^k x \dd x}{(x^2+a^2)(x^2+b^2)}. \]
By developing in power series around $\alpha=1$, we obtain
\[f(\alpha) \cos{\frac{\pi \alpha}{2}} = \frac{\pi}{2 (b^2-a^2)}
\sum_{n=0}^\infty \frac{ \log^n a - \log^n b}{n!} (\alpha -1)^n.\] By
differentiating $k$ times, \[ \sum_{j=0}^k \binom{k}{j} f^{(k-j)}(\alpha) \left(
\cos{\frac{\pi \alpha}{2}} \right)^{(j)} = \frac{\pi}{2 (b^2-a^2)}
\sum_{n=0}^\infty \frac{ \log^{n+k} a - \log^{n+k} b}{n!} (\alpha -1)^n.\] We
evaluate in $\alpha=1$, \[ \sum_{j=0\,\mathrm{(odd )} }^k  (-1)^{\frac{j+1}{2}}
\binom{k}{j} f^{(k-j)}(1) \left( \frac{\pi}{2} \right)^{j} = \frac{\pi (\log^k a
-  \log^k b)}{2 (b^2-a^2)}. \] As a consequence, we obtain \[f^{(k)}(1)=
\frac{1}{k+1} \sum_{j>1\,\mathrm{(odd )} }^{k+1}(-1)^{\frac{j+1}{2}}
\binom{k+1}{j} f^{(k+1-j)}(1) \left( \frac{\pi}{2} \right)^{j-1}+
\frac{\log^{k+1} a -  \log^{k+1} b}{ (k+1)(a^2-b^2)}. \] When $k=0$,
\[f^{(0)}(1)=f(1) = \frac{\log a -\log b}{a^2 -b^2} = \frac{\pi}{2}\frac{
P_{0}\left(\frac{2\log a}{\pi}\right) - P_{0}\left(\frac{2 \log
b}{\pi}\right)}{a^2-b^2}.\] The general result follows by induction on $k$ and
the definition of $P_k$. \qed

\section{Integrals and polylogarithms}

In order to understand how special values of zeta functions and L-series arise in our formulas, we are going to need the definition of polylogarithms, which can be found, for instance, in Goncharov's works, \cite{G1,G}:

\begin{defn} Multiple polylogarithms are defined as the power series
\[\Li_{n_1,\dots,n_m}(x_1,\dots, x_m) := \sum_{0 < k_1 < k_2 < \dots <k_m} \frac{x_1^{k_1}x_2^{k_2}\dots x_m^{k_m}}{k_1^{n_1}k_2^{n_2}\dots k_m^{n_m}}\]
which are convergent for $|x_i| < 1$. The length of a polylogarithm function is the number $m$ and its weight  is the number $w=n_1+ \dots + n_m$.
\end{defn}

\begin{defn} Hyperlogarithms are defined as the iterated integrals
\[\I_{n_1,\dots, n_m}(a_1:\dots:a_m:a_{m+1}) := \]
\[\int_0^{a_{m+1}} \underbrace{ \frac{\dd t}{t-a_1}\circ \frac{\dd t}{t} \circ \dots  \circ \frac{\dd t}{t}}_{n_1} \circ \underbrace{\frac{\dd t}{t-a_2} \circ \frac{\dd t}{t} \circ  \dots \circ \frac{\dd t}{t}}_{n_2} \circ \dots \circ \underbrace{\frac{\dd t}{t-a_m}\circ  \frac{\dd t}{t} \circ \dots \circ \frac{\dd t}{t}}_{n_m}\]
where $n_i$ are integers, $a_i$ are complex numbers, and
\[\int_0^{b_{k+1}} \frac{\dd t}{t-b_1} \circ \dots \circ \frac{\dd t}{t-b_k} =
\int_{0 \leq t_1 \leq \dots \leq t_k   \leq b_{k+1}} \frac {\dd t_1}{t_1-b_1}\,
\dots \, \frac{\dd t_k}{t_k - b_k}.\] \end{defn}

The value of the integral above only depends on the homotopy class of the path connecting $0$ and $a_{m+1}$ on $\Cset  \setminus \{ a_1, \dots, a_m\} $.

It is easy to see (for instance, in \cite{G}) that,
\begin{eqnarray*}
\I_{n_1,\dots, n_m}(a_1:\dots:a_m:a_{m+1}) & = & (-1)^m \Li_{n_1,\dots, n_m}\left(\frac{a_2}{a_1}, \frac{a_3}{a_2},\dots,\frac{a_m}{a_{m-1}}, \frac{a_{m+1}}{a_m}\right)\\\\
\Li_{n_1,\dots, n_m}(x_1, \dots, x_m) & = & (-1)^m \I_{n_1,\dots,  n_m}((x_1 \dots  x_m)^{-1}: \dots :x_m^{-1}: 1)
\end{eqnarray*}
which gives an analytic continuation of multiple polylogarithms.  Observe that
we recover the special value of the Riemann zeta function $\zeta(k)$ for $k \geq
2$ as $\Li_k(1)$, as well as $\Lf(\chi_{-4}, k) = - \frac{\ii}{2}( \Li_k(\ii) -
\Li_k(-\ii))$.

In order to express the results more clearly, we will establish some notation.

\begin{defn}
\begin{eqnarray*}
{\mathcal L}_r(\alpha) & := & \Li_r(\alpha) - \Li_r(-\alpha)\\\\
{\mathcal L}_{r,s}(\alpha,\alpha) & := & 2 ( \Li_{r,s}(\alpha, \alpha) - \Li_{r,s}(-\alpha, \alpha) + \Li_{r,s}(\alpha, -\alpha) - \Li_{r,s}(-\alpha, -\alpha) )
\end{eqnarray*}
\end{defn}
Note that the weight of any of the functions above is equal to the sum of its subindexes. This notation is the same as in \cite{L}.

Now we are ready to establish some technical results that will help us recognize special values of the Riemann zeta function and L-series out of integrals.
\begin{lem}\label{lemz} We have the following length-one identities:
\begin{eqnarray}
\int_0^1 \log^j x  \frac{\dd x}{x^2-1}& = & (-1)^{j+1} j! \left (1- \frac{1}{2^{j+1}} \right) \zeta(j+1)\\
\int_0^1 \log^j x  \frac{\dd x}{x^2+1}& = & \frac{(-1)^{j+1} j!}{2} \ii \mathcal{L}_{j+1}(\ii) = (-1)^j j! \Lf( \chi_{-4}, j+1)
\end{eqnarray}
and the following length-two  identities:
\begin{eqnarray}
\int_0^1  \int_0^x  \frac{\dd s}{s^2 -1} \circ \frac{\dd s}{s} \circ \frac{\dd s}{s} \log^{j} x \frac{ \dd x}{x^2 -1}  & = &  \frac{(-1)^j j!}{8} \mathcal{L}_{3,j+1}(1,1)\\
\int_0^1  \int_0^x  \frac{\dd s}{s^2 -1} \circ \frac{\dd s}{s} \circ \frac{\dd s}{s} \log^{j} x \frac{ \dd x}{x^2 + 1} & = & \frac{(-1)^{j+1} \ii  j!}{8} \mathcal{L}_{3,j+1}(\ii,\ii)
\end{eqnarray}
\end{lem}
\pf The idea is to translate the integral into hyperlogarithms. We use the fact that $\int_x^1 \frac{\dd s}{s} = - \log x$.
\[\int_0^1 \log^j x  \frac{\dd x}{x^2-1} = \frac{(-1)^j j!}{2} \int_0^1 \left( \frac{1}{x-1} - \frac{1}{x+1}\right) \dd x \circ \underbrace{\frac{\dd s}{ s} \circ \dots \circ \frac{\dd s}{s}}_{j\, \mathrm{times}} \]
The $j!$ occurs as a way to count the possible permutations of the variables $s$, since they are ordered in the hyperlogarithm integral.

\[= \frac{(-1)^{j+1} j!}{2} ( \Li_{j+1}(1) - \Li_{j+1}(-1)) =  \frac{(-1)^{j+1} j!}{2} 2 \left ( 1- \frac{1}{2^{j+1}} \right) \zeta(j+1)\]
The last equality is a consequence of the Euler product decomposition for the zeta function.
The second formula can be proved in a similar way.

Now, for the length-two identities, we do as before,
\[ \int_0^1  \int_0^x  \frac{\dd s}{s^2 -1} \circ \frac{\dd s}{s} \circ \frac{\dd s}{s} \log^{j} x \frac{ \dd x}{x^2-1 } \]
\[ =  \frac{ (-1)^j j!}{4}\int_0^1 \left( \frac{1}{s-1} - \frac{1}{s+1} \right) \dd s \circ \frac{\dd s}{s} \circ \frac{\dd s}{s} \circ \left( \frac{1}{x-1} - \frac{1}{x+1} \right) \dd x \circ \underbrace{\frac{\dd t}{t} \circ \dots \circ \frac{\dd t}{t}}_{j\, \mathrm{times}} \]
\[ =  \frac{(-1)^j j!}{4} ( \I_{3,j+1}(1:1:1) - \I_{3,j+1}(-1:1:1) + \I_{3,j+1}(-1:-1:1) - \I_{3,j+1}(1:-1:1) )\]
\[ =  \frac{(-1)^j j!}{4} ( \Li_{3,j+1}(1,1) - \Li_{3,j+1}(-1,1) + \Li_{3,j+1}(1,-1) - \Li_{3,j+1}(-1,-1) )\]
\[= \frac{ (-1)^j j!}{8} \mathcal{L}_{3,j+1}(1,1).\]
The other formula in the statement can be proved analogously. \qed

Now let us observe that the values $\mathcal{L}_{r,s}(1,1)$ for $r+s$ odd , can be expressed as combinations of values of $\zeta(k)$ for $2 \leq k \leq r+s$. This is possible because of the amazing formula (75) in \cite{BBB}, which claims:
\[ \Li_{r,s}(\rho,\sigma) = \frac{1}{2} \left( - \Li_{r+s}(\rho \sigma) + (1+(-1)^s) \Li_r(\rho) \Li_s(\sigma) \right)\]
\[ + \frac{(-1)^s}{2} \left( \binom{r+s-1}{r-1} \Li_{r+s}(\rho) + \binom{r+s-1}{s-1} \Li_{r+s} (\sigma) \right) \]
\[- \sum_{0 < k < \frac{r+s}{2} }\Li_{2k}(\rho \sigma)  (-1)^s\left(
\binom{r+s-2k-1}{r-1} \Li_{r+s-2k}(\rho) + \binom{r+s-2k-1}{s-1} \Li_{r+s-2k}
(\sigma) \right)\] \begin{equation}\label{BBB} \end{equation}
for $r+s$ odd , $\rho=\pm 1$, and $ \sigma=\pm 1$.

We will only need $\mathcal{L}_{r,s}(1,1)$ for $r=3$ and $s$ even.
\begin{prop} \label{propBBB} We have:
\[\mathcal{L}_{3,2h}(1,1) =\frac{2^{2h+3}-1}{2^{2h+1}} (h+1)(2h+1) \zeta(2h+3)  -\frac{2^{2h+1}-1}{2^{2h}} h(2h+5) \zeta(2) \zeta(2h+1) \]
\[- \sum_{k=2}^{h-1} \binom{2h-2k+2}{2} \frac{2^{2h-2k+3}-1}{2^{2h}} \zeta(2k) \zeta(2h-2k+3) \]
for $h \geq 2$, and
\[\mathcal{L}_{3,2}(1,1) = \frac{93}{4} \zeta(5) -\frac{21}{2} \zeta(2)
\zeta(3). \]

Expressing everything in terms of odd  special values of $\zeta$ and powers of $\pi$:
\[\mathcal{L}_{3,2h}(1,1) = \frac{2^{2h+3}-1}{2^{2h+1}} (h+1)(2h+1) \zeta(2h+3) -\frac{2^{2h+1}-1}{2^{2h}} h(2h+5) \frac{\pi^2}{6} \zeta(2h+1)\]
\[ - \sum_{k=2}^{h-1} \binom{2h-2k+2}{2} \frac{2^{2h-2k+3}-1}{2^{2h}} \frac{(-1)^{k-1}B_{2k} (2\pi)^{2k}}{2 (2k)!} \zeta(2h-2k+3)\]
for $h \geq 2$, and
\[\mathcal{L}_{3,2}(1,1) = \frac{93}{4} \zeta(5) -\frac{7}{4} \pi^2 \zeta(3). \]
\end{prop}
The proof of this statement is an easy application of formula (\ref{BBB}) together with the well-known formula
\[\zeta(2k) = \frac{(-1)^{k-1} B_{2k}(2 \pi)^{2k}}{2(2k)!}.\]

By applying formula (\ref{BBB}) we can also obtain the following result:
\begin{prop} \label{prop11} We have:
\[\int_0^1 \log(1+x) \log^j x \frac{\dd x}{x^2-1} = \frac{(-1)^j j!}{2} ( \Li_{1,
j+1} (-1, 1) -  \Li_{1, j+1} (1, -1) ) \]
and
\[  \Li_{1, 2h} (-1, 1) -  \Li_{1, 2h} (1, -1) = (2h -1) \frac{2^{2h+1} -1}{2^{2h+1}} \zeta(2h+1) - \log 2 \frac{2^{2h} -1}{2^{2h-1}} \zeta(2h) \]
\[- \sum_{k=1}^{h-1} \frac{(2^{2h+1-2k}-1)(2^{2k-1}-1)}{2^{2h-1}} \zeta(2k)
\zeta(2h+1-2k).\] In other words,
\[\int_0^1 \log(1+x) \log^{2h-1} x \frac{\dd x}{x^2-1}= -(2h-1)! (2h-1)\frac{2^{2h+1}-1}{2^{2h+2}} \zeta(2h+1)\]
\[+(-1)^{h-1}\frac{\log 2 (2^{2h}-1)  B_{2h}  \pi^{2h}}{4h}\]
\[ + \frac{(2h-1)!}{2}\sum_{k=1}^{h-1}
\frac{(2^{2h+1-2k}-1)(2^{2k-1}-1)}{2^{2h-2k}} \frac{(-1)^{k-1} B_{2k}}{(2k)!}
\pi^{2k} \zeta(2h+1-2k).\]
\end{prop}

\begin{prop} \label{log2} We have:
\[ \int_0^\infty \log(1+x^2) \log^{2h} x \frac{\dd x}{x^2+1}\]
\[= (-1)^h 2 E_{2h} \left( \frac{\pi}{2}\right)^{2h+1} \log 2\]
\[ + 2 \sum_{l=1}^{h} \frac{(2h)!}{(2h-2l)!} \left( 1- \frac{1}{2^{2l+1}}\right)
(-1)^{h-l} E_{2h-2l} \left( \frac{\pi}{2}\right)^{2h-2l+1} \zeta(2l+1)\]
where
$E_h$ is the $h$-th Euler number, $\frac{2\e^{x}}{\e^{2x}+1} = \sum_{n=0}^\infty
\frac{E_nx^n}{ n!}$.
\end{prop}
\pf The proof of this Proposition is again a
simple exercise in applying hyperlogarithms, \[\int_0^\infty \log(1+x^2)
\log^{2h} x \frac{\dd x}{x^2+1} = \int_0^\infty \int_0^1\frac{2x^2 t \dd t}{x^2t^2
+1}  \log^{2h} x \frac{\dd x}{x^2+1} \] \[ = 2 \int_0^1 \int_0^\infty \left(
\frac{1}{t^2x^2+1} - \frac{1}{x^2+1} \right)\log^{2h} x   \dd x \frac{t \dd
t}{1-t^2}.\]
Now we make the following change of variables: $y=tx$ in the first term but we
let $y=x$ in the second term, \[ = 2 \int_0^1 \int_0^\infty \frac{\dd y}{y^2+1}
(( \log y - \log t)^{2h}- t \log^{2h} y) \frac{ \dd t}{1-t^2} \] \[ =
2\int_0^\infty \frac{\log^{2h} y \dd y}{y^2 +1} \int_0^1 \frac{  \dd t}{1+ t}+ 2
\sum_{k=1}^{2h} \binom{2h}{k} (-1)^k \int_0^\infty \frac{\log^{2h-k} y \dd y}{y^2
+1} \int_0^1 \frac{  \log^k t \dd t}{1- t^2}.  \] We will apply Lemma \ref{lemz}.
But first, note that $\int_0^\infty \frac{\log^{n} y \dd y}{y^2 +1}=0$ for $n$
odd . Then we may write $k=2l$
\[ = 4 (2h)! \Lf(\chi_{-4}, 2h+1) \log 2 + 2 \sum_{l=1}^{h} \binom{2h}{2l}\]
\[  2 (2h-2l)! \Lf(\chi_{-4}, 2h-2l+1) (2l)!\left( 1- \frac{1}{2^{2l+1}}\right) \zeta(2l+1) \]
\[ = 4 (2h)! \Lf(\chi_{-4}, 2h+1) \log 2 + 4 \sum_{l=1}^{h} (2h)! \left( 1-
\frac{1}{2^{2l+1}}\right)  \Lf(\chi_{-4}, 2h-2l+1)   \zeta(2l+1). \] The proof
of the statement is now just an application of the well-known formula \[
\Lf(\chi_{-4}, 2k+1) = \frac{(-1)^k E_{2k} \left( \frac{\pi}{2}\right)^{2k+1}}{2
(2k)!}.\] \qed

\begin{prop} \label{li2} We have:
\[ \frac{1}{2\ii} \int_{0}^\infty  (\left ( \Li_2\left ( \ii x \right ) - \Li_2(-\ii x)\right ) \log^{2h} x  \frac{\dd x}{x^2+1}\]
\[ =\sum_{l=0}^{h} B_{2l} \frac{(2h)!}{(2l)!} (2^{2l-1}-1) (-1)^{l+1} \pi^{2l}
(h-l+1) \left( \frac{2^{2h+3-2l} - 1}{2^{2h+1}} \right)\zeta(2h+3-2l).\]
\end{prop}
\pf Note that
\[ \frac{1}{2\ii} \int_{0}^\infty  (\left ( \Li_2\left ( \ii x \right ) - \Li_2(-\ii x)\right ) \log^{2h} x  \frac{\dd x}{x^2+1} \]
\[= \frac{1}{2\ii} \int_{0}^\infty  \int_0^1 \frac{2 \ii x \dd t}{t^2x^2+1} \circ
\frac{\dd t}{t}  \log^{2h} x  \frac{\dd x}{x^2+1}  = -\int_{0}^1  \int_0^\infty
\frac{x \log^{2h} x  \dd x  }{(t^2x^2+1)(x^2+1)} \log t \dd t. \]
We are now ready to apply Proposition \ref{prop} \[= -\int_{0}^1  \left(
\frac{\pi}{2} \right)^{2h+1} \frac{P_{2h} \left ( -\frac{2 \log t}{\pi} \right)
}{1 - t^2} \log t \dd t . \] Applying Proposition \ref{propP} from the Appendix,
\[= \int_{0}^1  \left( \frac{\pi}{2} \right)^{2h+1} \frac{2}{2h+1}
\sum_{l=0}^{h} B_{2l} \binom{2h+1}{2l} (2^{2l-1}-1) (-1)^l  \frac{\left (
-\frac{2 \log t}{\pi} \right)^{2h+1-2l }}{1 - t^2} \log t \dd t  \] \[=
\int_{0}^1  \frac{2}{2h+1} \sum_{l=0}^{h} B_{2l} \binom{2h+1}{2l} (2^{2l-1}-1)
(-1)^{l+1} \left( \frac{\pi}{2} \right )^{2l} \frac{ \log ^{2h+2-2l } t }{1 -
t^2} \dd t. \] Then, apply Lemma \ref{lemz} \[=  2 \sum_{l=0}^{h} B_{2l}
\frac{(2h)!}{(2l)!} (2^{2l-1}-1) (-1)^{l+1} \left( \frac{\pi}{2} \right )^{2l}
(2h-2l+2) \left( 1 - \frac{1}{2^{2h+3-2l}} \right)\zeta(2h+3-2l)  \] \[=
\sum_{l=0}^{h} B_{2l} \frac{(2h)!}{(2l)!} (2^{2l-1}-1) (-1)^{l+1} \pi^{2l}
(h-l+1) \left( \frac{2^{2h+3-2l} - 1}{2^{2h+1}} \right)\zeta(2h+3-2l).  \] \qed

\section{An identity for symmetrical polynomials}
For dealing with the polynomials $P_k$, we will need to manage certain identities of symmetric polynomials. More specifically, we are going to use the following result:
\begin{prop}\label{propsym}
\begin{eqnarray*}
2n (-1)^l s_{n-l}(2^2,\dots,(2n-2)^2) & = & \sum_{h=l}^n (-1)^{h}  \binom{2h}{2l-1} s_{n-h}(1^2, \dots , (2n-1)^2) \\
(2n+1) (-1)^l s_{n-l}(1^2, \dots , (2n-1)^2) & = & \sum_{h=l}^{n} (-1)^{h} \binom{2h+1}{2l} s_{n-h}(2^2, \dots , (2n)^2)
\end{eqnarray*}
\end{prop}

\pf These equalities are easier to prove if we think of the symmetric functions
as coefficients of certain polynomials, as in equation (\ref{eq:sym}).

In order to prove the first equality, multiply by $x^{2l}$ on both sides and add
for $l=1, \dots, n$:
\[2n \sum_{l=1}^n s_{n-l}(2^2,\dots,(2n-2)^2) (-1)^l x^{2l} \]
\[=\sum_{l=1}^n\sum_{h=l}^n (-1)^{h}  \binom{2h}{2l-1} s_{n-h}(1^2, \dots ,
(2n-1)^2) x^{2l}.\]

The statement we have to prove becomes:
\[2n \prod_{j=0}^{n-1} ((2j)^2-x^2) = \sum_{h=1}^n (-1)^h s_{n-h}(1^2, \dots ,
(2n-1)^2) \sum_{l=1}^h   \binom{2h}{2l-1} x^{2l}.
\]\begin{equation}\label{eq:sym1} \end{equation}

The right side of (\ref{eq:sym1}) is
\[= \sum_{h=0}^n (-1)^h s_{n-h}(1^2, \dots , (2n-1)^2) \frac{x}{2} ((x+1)^{2h} -
(x-1)^{2h}) \]
\[=\frac{x}{2} \left(\prod_{j=1}^n((2j-1)^2 -(x+1)^2) -\prod_{j=1}^n((2j-1)^2 -(x-1)^2) \right) \]
\[=\frac{x}{2} \left(\prod_{j=1}^n(2j+x)(2j-2-x) -\prod_{j=1}^n(2j-x)(2j-2 +x) \right) \]
\[= ((-x)(2n+x) - x(2n-x)) \frac{x}{2} \prod_{j=1}^{n-1} ((2j)^2 - x^2)= 2n \prod_{j=0}^{n-1} ((2j)^2-x^2)\]
so equation (\ref{eq:sym1}) is true.

In order to prove the second equality, we apply a similar process. First
multiply by $x^{2l+1}$ on both sides and add  for $l=1, \dots, n$:
\[(2n+1) \sum_{l=1}^n s_{n-l}(1^2,\dots,(2n-1)^2) (-1)^l x^{2l+1} \]
\[=\sum_{l=1}^n\sum_{h=l}^n (-1)^{h}  \binom{2h+1}{2l} s_{n-h}(2^2, \dots , (2n)^2)
x^{2l+1}.\]

Hence, we have to prove:
\[(2n+1) x \prod_{j=1}^{n} ((2j-1)^2-x^2) = \sum_{h=1}^n (-1)^h s_{n-h}(2^2,
\dots , (2n)^2) \sum_{l=1}^h   \binom{2h+1}{2l} x^{2l+1}.\] \begin{equation}
\label{eq:sym2} \end{equation}

The right side is
\[= \sum_{h=0}^n (-1)^h s_{n-h}(2^2, \dots , (2n)^2) \frac{x}{2} ((x+1)^{2h+1} - (x-1)^{2h+1}) \]
\[=\frac{x}{2} \left((x+1) \prod_{j=1}^n((2j)^2 -(x+1)^2) -(x-1)\prod_{j=1}^n((2j)^2 -(x-1)^2) \right) \]
\[=\frac{x}{2} \left((x+1)\prod_{j=1}^n(2j+1 +x)(2j-1-x) -(x-1)\prod_{j=1}^n(2j-1 +x)(2j+1 -x) \right) \]
\[= ((2n+1+x) + (2n+1-x)) \frac{x}{2} \prod_{j=1}^{n} ((2j-1)^2 - x^2)\]
\[= (2n+1) x\prod_{j=1}^{n} ((2j-1)^2-x^2)\] thus proving equation
(\ref{eq:sym2}). \qed

\section{Description of the general method}
We will prove our main result by first examining a general situation and  then
specializing to the particular families of the statement.

Let $P_\alpha \in \Cset [{\bf x}]$ such that its coefficients depend polynomially on a parameter $\alpha \in \Cset$. We replace $\alpha$ by $\left(\frac{x_1-1}{x_1+1}\right) \dots \left(\frac{x_n-1}{x_n+1}\right)$ and obtain a new polynomial $\tilde{P} \in \Cset[{\bf x}, x_1, \dots, x_n]$. By definition of Mahler measure, it is easy to see that
\[m(\tilde{P}) = \frac{1}{(2 \pi \ii)^n} \int_{\TT^n}
m\left(P_{\left(\frac{x_1-1}{x_1+1}\right) \dots
\left(\frac{x_n-1}{x_n+1}\right)}\right) \frac{\dd x_1}{x_1} \dots \frac{\dd
x_n}{x_n}. \] We perform a change of variables to polar coordinates, $x_j =
\e^{\ii \theta_j}$: \[= \frac{1}{(2 \pi )^n} \int_{-\pi}^{\pi} \dots
\int_{-\pi}^{\pi}  m\left(P_{\ii^n \tan \left(\frac{\theta_1}{2}\right) \dots
\tan \left(\frac{\theta_n}{2}\right) }\right) \dd \theta_1 \dots \dd \theta_n. \]
Set $x_i = \tan \left( \frac{\theta_i}{2} \right) $. We get, \[= \frac{1}{\pi
^n} \int_{-\infty}^{\infty} \dots \int_{-\infty}^{\infty} m\left(P_{\ii^n x_1
\dots x_n }\right) \frac{\dd x_1}{x_1^2+1}\dots \frac{\dd x_n}{x_n^2+1}\]
\[ =\frac{2^n}{\pi ^n} \int_0^{\infty} \dots \int_0^{\infty} m\left(P_{\ii^n x_1
\dots x_n }\right) \frac{\dd x_1}{x_1^2+1}\dots \frac{\dd x_n}{x_n^2+1} . \]
Making one more change, $\hat{x}_1 = x_1, \dots, \hat{x}_{n-1} = x_1 \dots
x_{n-1}, \hat{x}_n = x_1 \dots x_n $: \[= \frac{2^n}{\pi ^n} \int_0^{\infty}
\dots \int_0^{\infty}  m\left(P_{\ii^n \hat{x}_n }\right) \frac{\hat{x}_1 \dd
\hat{x}_1}{\hat{x}_1^2+1}\frac{\hat{x}_2 \dd
\hat{x}_2}{\hat{x}_2^2+\hat{x}_1^2}\dots \frac{\hat{x}_{n-1}\dd
\hat{x}_{n-1}}{\hat{x}_{n-1}^2+\hat{x}_{n-2}^2} \frac{\dd
\hat{x}_n}{\hat{x}_n^2+\hat{x}_{n-1}^2}. \]

We need to compute this integral. In most of our cases, the Mahler measure of
$P_\alpha$ depends only on the absolute value of $\alpha$. If not, for each
$n$ we may modify $P$, such that it absorbs the number  $\ii^n$. From now on, we
will write $m(P_x)$ instead of $m(P_{\ii^n x})$ to simplify notation.

By iterating Proposition \ref{prop}, the above integral can be written as a
linear combination, with coefficients that are rational numbers and powers of
$\pi$ in such a way that the weights are homogeneous,  of integrals of the form
\[ \int_0^\infty  m\left(P_{x}\right)\log^j x \frac{ \dd x}{x^2 \pm 1}.  \]

It is easy to see that $j$ is even iff $n$ is odd  and the corresponding sign in that case is "+".

We are going to compute these coefficients.

Let us establish some convenient notation:

\begin{defn} Let $a_{n,h} \in \Qset$ be defined for $ n \geq 1$ and $h=0, \dots, n-1$ by
 \[ \int_0^{\infty} \dots \int_0^{\infty}  m\left(P_{x_1 }\right) \frac{x_{2n} \dd x_{2n}}{x_{2n}^2+1}\frac{x_{2n-1} \dd x_{2n-1}}{x_{2n-1}^2+ x_{2n}^2}\dots \frac{\dd x_1}{x_1^2+ x_ 2^2} \]
\begin{equation}
= \sum_{h=1}^{n} a_{n,h-1} \left(\frac{\pi}{2}\right)^{2n-2h}\int_0^\infty
m\left(P_{x}\right)\log^{2h-1} x \frac{ \dd x}{x^2 - 1}. \end{equation}

Let $b_{n,h} \in \Qset$ be defined for $ n \geq 0$ and $h=0, \dots, n$ by
\[ \int_0^{\infty} \dots \int_0^{\infty}  m\left(P_{x_1 }\right) \frac{x_{2n+1} \dd x_{2n+1}}{x_{2n+1}^2+1}\frac{x_{2n} \dd x_{2n}}{x_{2n}^2+ x_{2n+1}^2}\dots \frac{\dd x_1}{x_1^2+ x_ 2^2} \]
\begin{equation}
= \sum_{h=0}^{n} b_{n,h} \left(\frac{\pi}{2}\right)^{2n-2h}\int_0^\infty
m\left(P_{x}\right)\log^{2h} x \frac{ \dd x}{x^2 + 1}. \end{equation}
\end{defn}

We claim:

\begin{lem} \label{lemma13}
\begin{eqnarray}
\sum_{h=0}^n b_{n,h} x^{2h} & =  & \sum_{h=1}^{n} a_{n,h-1} \left(P_{2h-1}\left(x\right) - P_{2h-1}\left( \ii \right)\right) \label{eq:idab1}\\
\sum_{h=1}^{n+1} a_{n+1,h-1} x^{2h-1} & = &  \sum_{h=0}^{n} b_{n,h} P_{2h}\left(x\right)   \label{eq:idab2}
\end{eqnarray}
\end{lem}

\pf First observe that
\[  \sum_{h=0}^{n} b_{n,h} \left(\frac{\pi}{2}\right)^{2n-2h}\int_0^\infty  m\left(P_{x}\right)\log^{2h} x \frac{ \dd x}{x^2 + 1} \]
\begin{equation} \label{eq:ab1}
=  \sum_{h=1}^{n} a_{n,h-1} \left(\frac{\pi}{2}\right)^{2n-2h}\int_0^\infty
\int_0^\infty  m\left(P_{x}\right)  y \log^{2h-1} y \frac{ \dd y}{y^2 - 1}
\frac{\dd x}{x^2+y^2}. \end{equation}

But
\[ \int_0^\infty \frac{y \log^{2h-1} y  \dd y}{(y^2+x^2)(y^2 - 1)}  =  \left(\frac{\pi}{2}\right)^{2h}\frac{ P_{2h-1}\left(\frac{2\log x}{\pi}\right) - P_{2h-1}\left( \ii \right)}{x^2 + 1}\]
by applying Proposition \ref{prop} for $a = x$ and $b= \ii$.

The right side of equation (\ref{eq:ab1}) becomes
\[  =  \sum_{h=1}^{n} a_{n,h-1} \left(\frac{\pi}{2}\right)^{2n} \int_0^\infty
m\left(P_{x}\right)  \left(P_{2h-1}\left(\frac{2\log x}{\pi}\right) -
P_{2h-1}\left( \ii \right)\right) \frac{\dd x} {x^2 + 1}.\]

As a consequence, equation  (\ref{eq:ab1})  translates into the polynomial identity (\ref {eq:idab1}).

On the other hand,
\[  \sum_{h=1}^{n+1} a_{n+1,h-1} \left(\frac{\pi}{2}\right)^{2n+2-2h}\int_0^\infty  m\left(P_{x}\right)\log^{2h-1} x \frac{ \dd x}{x^2 - 1} \]
\begin{equation}\label{eq:ab2}
 =  \sum_{h=0}^{n} b_{n,h} \left(\frac{\pi}{2}\right)^{2n-2h}\int_0^\infty
\int_0^\infty  m\left(P_{x}\right)  y \log^{2h} y \frac{ \dd y}{y^2 + 1}
\frac{\dd x}{x^2+y^2}. \end{equation}

But
\[ \int_0^\infty \frac{y \log^{2h} y  \dd y}{(y^2+x^2)(y^2 + 1)}  =  \left(\frac{\pi}{2}\right)^{2h+1}\frac{ P_{2h}\left(\frac{2\log x}{\pi}\right) - P_{2h}\left( 0 \right)}{x^2 - 1}\]
by applying Proposition \ref{prop} for $a = x$ and $b= 1$.

So the right side of (\ref{eq:ab2}) becomes
\[  =  \sum_{h=0}^{n} b_{n,h} \left(\frac{\pi}{2}\right)^{2n+1} \int_0^\infty  m\left(P_{y}\right)  P_{2h}\left(\frac{2\log x}{\pi}\right)  \frac{\dd x} {x^2 - 1}\]
which translates into the identity (\ref{eq:idab2}). \qed

\begin{thm} \label{thm14} We have:
\begin{equation}
\sum_{h=0}^{n-1} a_{n,h} x^{2h}= \frac{(x^2+2^2) \dots (x^2+(2n-2)^2)}{(2n-1)!}
\end{equation}
for $ n \geq 1$ and $h=0, \dots, n-1$, and
\begin{equation}
\sum_{h=0}^n b_{n,h} x^{2h}= \frac{(x^2+1^2) \dots (x^2+(2n-1)^2)}{(2n)!}
\end{equation}
for $ n \geq 0$ and $h=0, \dots, n$.

In other words,
\begin{eqnarray}
a_{n,h} & = &  \frac{s_{n-1-h}(2^2, \dots, (2n-2)^2)}{(2n-1)!}\\
b_{n,h} & = & \frac{s_{n-h}(1^2, \dots, (2n-1)^2)}{(2n)!}
\end{eqnarray}
\end{thm}

\pf For $2n+1=1$, $n=0$ and the integral becomes
\[\int_0^{\infty} m\left(P_{x}\right) \frac{\dd x}{x^2+1}\]
so $b_{0,0}=1$.

For $2n=2$, $n=1$ and we have
\[\int_0^{\infty} \int_0^{\infty} m\left(P_{x}\right) \frac{y \dd y}{y^2+1}\frac{\dd x}{x^2+y^2}=
\int_0^{\infty} m\left(P_{x}\right) \frac{\log x \dd x}{x^2-1}\]
so $a_{1,0} = 1$.

Then the statement is true for the first two cases.

We proceed by induction. Suppose that
\[ a_{n,h} = \frac{s_{n-1-h}(2^2, \dots, (2n-2)^2)}{(2n-1)!}.\]
We have  to prove that
\[ b_{n,h} =  \frac{s_{n-h}(1^2, \dots, (2n-1)^2)}{(2n)!}.\]

By Lemma \ref{lemma13}, it is enough to prove that
\[\sum_{h=0}^n s_{n-h}(1^2, \dots, (2n-1)^2)  x^{2h}=   2n \sum_{h=1}^{n}  s_{n
-h}(2^2, \dots, (2n-2)^2) \left(P_{2h-1}\left(x\right) - P_{2h-1}\left( \ii
\right)\right).\]
\begin{equation}\label{eq:inducab}
\end{equation}

Recall equation (\ref{defiP}) that defines the polynomials $P_k$, from which the following identity may be deduced:
\begin{equation} \label{eqeven}
 x^{2h} =  \sum_{k=0}^{h-1}(-1)^{k} \binom{2h}{2k+1} P_{2h - 2k-1}(x).
\end{equation}

Multiplying equation (\ref{eqeven}) by $ s_{n-h}(1^2, \dots , (2n-1)^2)   $ and add ing, we get
\[\sum_{h=0}^n s_{n-h}(1^2, \dots , (2n-1)^2) x^{2h} \]
\[ = \sum_{h=1}^n s_{n-h}(1^2, \dots , (2n-1)^2)  \sum_{k=0}^{h-1}(-1)^k
\binom{2h}{2k+1} P_{2h - 2k-1}(x) + s_n(1^2, \dots , (2n-1)^2).\]

Now let us evaluate the above equality at $x=\ii$, we obtain
\[\sum_{h=0}^n s_{n-h}(1^2, \dots , (2n-1)^2) (-1)^h\]
\[  = \sum_{h=1}^n s_{n-h}(1^2, \dots , (2n-1)^2)  \sum_{k=0}^{h-1}(-1)^k
\binom{2h}{2k+1} P_{2h - 2k-1}(\ii) + s_n(1^2, \dots , (2n-1)^2).\] But
\[\sum_{h=0}^n s_{n-h}(1^2, \dots , (2n-1)^2) (-1)^h =(x +1^2) \dots (x
+(2n-1)^2)|_{x=-1} =0, \] from where
\[\sum_{h=0}^n s_{n-h}(1^2, \dots , (2n-1)^2) x^{2h} \]
\[ = \sum_{h=1}^n s_{n-h}(1^2, \dots , (2n-1)^2)  \sum_{k=0}^{h-1}(-1)^k
\binom{2h}{2k+1} (P_{2h - 2k-1}(x) -P_{2h - 2k-1}(\ii)) . \]

Let $l =h-k$, then this becomes
\[ = \sum_{h=1}^n s_{n-h}(1^2, \dots , (2n-1)^2)  \sum_{l=1}^{h}(-1)^{h-l} \binom{2h}{2l-1} (P_{2l -1}(x) -P_{2l-1}(\ii))  \]
\[ = \sum_{l=1}^n \left(\sum_{h=l}^n (-1)^h  \binom{2h}{2l-1} s_{n-h}(1^2, \dots
, (2n-1)^2) \right) (-1)^l(P_{2l -1}(x) -P_{2l-1}(\ii)) , \] and equality
(\ref{eq:inducab}) is proved by applying Proposition \ref{propsym}.

Now suppose that
\[ b_{n,h} = \frac{s_{n-h}(1^2, \dots, (2n-1)^2)}{(2n)!},\]
we want to see that
\[ a_{n+1,h} =  \frac{s_{n-h}(2^2, \dots, (2n)^2)}{(2n+1)!}.\]

Then it is enough to prove that
\[\sum_{h=0}^{n} s_{n-h}(2^2, \dots, (2n)^2)  x^{2h+1}=   (2n+1) \sum_{h=0}^{n}
s_{n -h}(1^2, \dots, (2n-1)^2) P_{2h}\left(x\right)\]
\begin{equation}\label{eq:inducba} \end{equation}
by Lemma \ref{lemma13}.

Equation (\ref{defiP}) implies
\[ x^{2h+1} =  \sum_{k=0}^{h}(-1)^k \binom{2h+1}{2k+1} P_{2h - 2k}(x), \]
and so,
\[\sum_{h=0}^{n} s_{n-h}(2^2, \dots, (2n)^2) x^{2h+1} \]
\[ = \sum_{h=0}^{n}
s_{n-h}(2^2, \dots , (2n)^2)  \sum_{k=0}^{h}(-1)^k \binom{2h+1}{2k+1} P_{2h -
2k}(x). \]

Let $l =h-k$, then
\[= \sum_{h=0}^{n} s_{n-h}(2^2, \dots , (2n)^2)  \sum_{l=0}^{h}(-1)^{h-l} \binom{2h+1}{2l} P_{2l}(x) \]

\[= \sum_{l=0}^{n}\left( \sum_{h=l}^{n} (-1)^h \binom{2h+1}{2l} s_{n-h}(2^2, \dots , (2n)^2)  \right) (-1)^l P_{2l}(x) \]
which proves (\ref{eq:inducba}) by Proposition \ref{propsym}. \qed

\section{Proof of the main theorem }

In the last section we managed to express the Mahler measure of $\tilde{P}$ as a linear combination of functions that depend on the Mahler measure of $P_\alpha$.  We are now ready to apply that machinery to the specific families of polynomials.
At this point we need to strongly use the formulas for the Mahler measure of each particular polynomial $P_\alpha$.

\bigskip
(i) $P_\alpha(z) =1+\alpha z$.
\[m(1+\alpha z) = \log^+|\alpha|\]
This is the simplest possible case. For the even case we get
\[ \pi^{2n} m\left( 1+  \left( \frac{1-x_1}{1+x_1}\right) \dots \left(
\frac{1-x_{2n}}{1+x_{2n}}\right) z \right)\]
\[= 2^{2n} \sum_{h=1}^{n} a_{n,h-1} \left(\frac{\pi}{2}\right)^{2n-2h}\int_0^\infty  \log^+ x \log^{2h-1} x \frac{ \dd x}{x^2 - 1} \]
\[= \sum_{h=1}^n  \frac{s_{n-h}(2^2, \dots, (2n-2)^2)}{(2n-1)!} 2^{2h}
\pi^{2n-2h} \int_1^\infty \log^{2h} x \frac{\dd x}{x^2-1}.\] Now set
$y=\frac{1}{x}$, \[= \sum_{h=1}^n  \frac{s_{n-h}(2^2, \dots, (2n-2)^2)}{(2n-1)!}
2^{2h} \pi^{2n-2h} \int_0^1 \log^{2h} y \frac{\dd y}{1-y^2}.\] If we apply  Lemma
\ref{lemz}, we obtain \[=\sum_{h=1}^n  \frac{s_{n-h}(2^2, \dots,
(2n-2)^2)}{(2n-1)!} 2^{2h} \pi^{2n-2h}(2h)!\left( 1- \frac{1}{2^{2h+1}} \right)
\zeta(2h+1)\] \[=\sum_{h=1}^n  \frac{s_{n-h}(2^2, \dots, (2n-2)^2)}{(2n-1)!}
\frac{(2h)!(2^{2h+1}-1)}{2}  \pi^{2n-2h} \zeta(2h+1).\]

For the odd  case we get
\[ \pi^{2n+1} m\left( 1+  \left( \frac{1-x_1}{1+x_1}\right) \dots \left(
\frac{1-x_{2n+1}}{1+x_{2n+1}}\right) z \right)\]
\[=2^{2n+1} \sum_{h=0}^{n} b_{n,h} \left(\frac{\pi}{2}\right)^{2n-2h}\int_0^\infty  \log^+ x \log^{2h} x \frac{ \dd x}{x^2 + 1} \]
\[= \sum_{h=0}^n  \frac{s_{n-h}(1^2, \dots, (2n-1)^2)}{(2n)!} 2^{2h+1}
\pi^{2n-2h} \int_1^\infty \log^{2h+1} x \frac{\dd x}{x^2+1}.\] Now set
$y=\frac{1}{x}$, \[= - \sum_{h=0}^n  \frac{s_{n-h}(1^2, \dots, (2n-1)^2)}{(2n)!}
2^{2h+1} \pi^{2n-2h} \int_0^1 \log^{2h+1} y \frac{\dd y}{y^2+1}.\] Now apply
Lemma \ref{lemz}, \[=  \sum_{h=0}^n  \frac{s_{n-h}(1^2, \dots, (2n-1)^2)}{(2n)!}
(2h+1)! 2^{2h+1} \pi^{2n-2h} \Lf(\chi_{-4}, 2h+2).\]

\bigskip

(ii) $P_\alpha(x,y,z) =(1+x)+\alpha (1+y)z$.
\bigskip

This Mahler measure was computed by Smyth \cite{B1,S2},
\[\pi^2 m((1+x)+\alpha (1+y)z) = \left \{\begin{array}{ll} 2 \mathcal{L}_3\left(|\alpha| \right) & \mathrm{for}\quad |\alpha| \leq 1\\ \\
\pi^2 \log |\alpha| + 2 \mathcal{L}_3\left(|\alpha|^{-1} \right) & \mathrm{for}\quad |\alpha| > 1
\end{array} \right .
\]
Where
\[\mathcal{L}_3\left(\alpha \right) = - \frac{2}{\alpha}  \int_0^1 \frac{\dd
t}{t^2 - \frac{1}{\alpha^2}}\circ \frac{\dd t}{t} \circ \frac{\dd t}{t}. \] Now
set $s = t \alpha$, \[=-2 \int_0^\alpha \frac{\dd s}{s^2 -1} \circ \frac{\dd s}{s}
\circ \frac{\dd s}{s}. \]
We obtain
\[ \pi^{2n+2} m\left(  1+x +  \left( \frac{1-x_1}{1+x_1}\right) \dots \left(
\frac{1-x_{2n}}{1+x_{2n}}\right)(1+y) z \right)\]
\[= 2^{2n} \sum_{h=1}^{n} a_{n,h-1} \left(\frac{\pi}{2}\right)^{2n-2h}\left(-4 \int_0^1 \int_0^x  \frac{\dd s}{s^2 -1} \circ \frac{\dd s}{s} \circ \frac{\dd s}{s} \log^{2h-1} x \frac{ \dd x}{x^2 - 1} \right . \]
\[\left . +\pi^2 \int_1^\infty  \log^{2h} x \frac{ \dd x}{x^2 - 1} -4   \int_1^\infty    \int_0^\frac{1}{x} \frac{\dd s}{s^2 -1} \circ \frac{\dd s}{s} \circ \frac{\dd s}{s} \log^{2h-1} x \frac{ \dd x}{x^2 - 1} \right)\]
\[=  \sum_{h=1}^{n} \frac{s_{n-h}(2^2, \dots, (2n-2)^2)}{(2n-1)!} 2^{2h} \pi ^{2n-2h} \]
\[ \left( -8 \int_0^1  \int_0^x  \frac{\dd s}{s^2 -1} \circ \frac{\dd s}{s} \circ
\frac{\dd s}{s} \log^{2h-1} x \frac{ \dd x}{x^2 - 1} - \pi^2\int_0^1  \log^{2h} x
\frac{ \dd x}{x^2 - 1} \right).\]
By Lemma \ref{lemz}, \[=  \sum_{h=1}^{n}
\frac{s_{n-h}(2^2, \dots, (2n-2)^2)}{(2n-1)!} 2^{2h} \pi ^{2n-2h} \] \[
\left((2h-1)! \mathcal{L}_{3,2h}(1,1) + (2h)! \left( 1 -
\frac{1}{2^{2h+1}}\right)\pi^2 \zeta(2h+1) \right).\]

But if we apply Proposition \ref{propBBB}, we get just combinations of the Riemann zeta function:
\[=  \sum_{h=1}^{n} \frac{s_{n-h}(2^2, \dots, (2n-2)^2)}{(2n-1)!} \pi ^{2n-2h} \]
\[(2h-1)!\sum_{k=0}^{h-1} \binom{2h-2k+2}{2} (2^{2h-2k+3}-1) \frac{(-1)^{k}
B_{2k}(2\pi)^{2k}}{2(2k)!} \zeta(2h-2k+3). \]

Now set $t=h-k$ and change the order of the sums:
\[=  \sum_{t=1}^{n}\frac{(2t+2)!(2^{2t+3}-1)}{8}\]
\[  \left(\sum_{k=0}^{n-t} \frac{s_{n-t-k}(2^2, \dots, (2n-2)^2)}{(2n-1)!}
\binom{2t+2k}{2t}  (-1)^{k} \frac{2^{2k} }{t+k} B_{2k}\right)\pi
^{2n-2t}\zeta(2t+3).\]

The odd  case is
\[ \pi^{2n+3} m\left(  1+x +  \left( \frac{1-x_1}{1+x_1}\right) \dots \left(
\frac{1-x_{2n+1}}{1+x_{2n+1}}\right)(1+y) z \right)\]
\[= 2^{2n+1} \sum_{h=0}^{n} b_{n,h} \left(\frac{\pi}{2}\right)^{2n-2h}\left(-4 \int_0^1 \int_0^x  \frac{\dd s}{s^2 -1} \circ \frac{\dd s}{s} \circ \frac{\dd s}{s} \log^{2h} x \frac{ \dd x}{x^2 +1} \right . \]
\[\left . +\pi^2 \int_1^\infty  \log^{2h+1} x \frac{ \dd x}{x^2 + 1} -4   \int_1^\infty    \int_0^\frac{1}{x} \frac{\dd s}{s^2 -1} \circ \frac{\dd s}{s} \circ \frac{\dd s}{s} \log^{2h} x \frac{ \dd x}{x^2 + 1} \right)\]
\[=  \sum_{h=0}^{n} \frac{s_{n-h}(1^2, \dots, (2n-1)^2)}{(2n)!} 2^{2h+1} \pi ^{2n-2h} \]
\[ \left( -8 \int_0^1  \int_0^x  \frac{\dd s}{s^2 -1} \circ \frac{\dd s}{s} \circ
\frac{\dd s}{s} \log^{2h} x \frac{ \dd x}{x^2 + 1} - \pi^2 \int_0^1  \log^{2h+1} x
\frac{ \dd x}{x^2 + 1} \right).\] By Lemma \ref{lemz}, \[=  \sum_{h=0}^{n}
\frac{s_{n-h}(1^2, \dots, (2n-1)^2)}{(2n)!} 2^{2h+1} \pi ^{2n-2h} \] \[
\left(\ii (2h)! \mathcal{L}_{3,2h+1}(\ii,\ii) + (2h+1)!\pi^2 \Lf(\chi_{-4},
2h+2)\right).\]

We should observe that it would be nice to have a simpler expression for $\mathcal{L}_{3,2h+1}(\ii,\ii)$. In fact, we believe that this number should be somehow related to $\Lf(\chi_{-4}, k)$ (in a result analogous to equation (\ref{BBB})), but we have been unable to find such a relation.

\bigskip

(iii) $P_\alpha(z) =1+\alpha x + (1-\alpha) y $.
\bigskip

This Mahler measure is a particular case of an example computed by Cassaigne
and Maillot \cite{M}. This case is different from the previously studied cases
due to the fact that the Mahler measure of this polynomial does not just depend
on the absolute value of the parameter $\alpha$, it also depends on the argument
of $\alpha$. This fact makes the application of the general method a little bit
more subtle. We will use \[\pi m(1+\alpha x +(1-\alpha) y  )\]
\[ = | \arg \alpha |
\log |1-\alpha| + | \arg (1- \alpha)|  \log |\alpha| + \left \{
\begin{array}{cc} D ( \alpha) & \mbox{if} \, \im (\alpha) \geq 0\\ \\  D (
\bar{\alpha}) & \mbox{if} \, \im (\alpha) < 0 \end{array} \right . \] The
deduction of this formula can be found in \cite{L}. For the even case we need to
use the formula for the case in which the parameter $\alpha$ is real, \[
m(1+\alpha x +(1-\alpha) y  ) = \left \{ \begin{array}{cc}  \log^+ \alpha   &
\mbox{if} \, \alpha > 0 \\\\ \log (1 - \alpha) &  \mbox{if} \, \alpha < 0
\end{array} \right .\] Then \[ \pi^{2n + 1} m\left( 1+  \left(
\frac{1-x_1}{1+x_1}\right) \dots \left( \frac{1-x_{2n}}{1+x_{2n}}\right) x +
\left(1 -  \left( \frac{1-x_1}{1+x_1}\right) \dots \left(
\frac{1-x_{2n}}{1+x_{2n}}\right) \right) y  \right)\] \[= \pi 2^{2n}
\sum_{h=1}^{n} a_{n,h-1}
\left(\frac{\pi}{2}\right)^{2n-2h}\frac{1}{2}\int_{-\infty}^\infty m(P_{(-1)^n
x}) \log^{2h-1}| x| \frac{ \dd x}{x^2 - 1}.\] Note that we have taken into
account that the formula depends on the sign of the parameter. \[= \sum_{h=1}^n
\frac{s_{n-h}(2^2, \dots, (2n-2)^2)}{(2n-1)!} 2^{2h} \pi^{2n+1-2h} \int_0^\infty
\frac{1}{2}(\log^+ x + \log (1+x)) \log^{2h-1} x \frac{ \dd x}{x^2 - 1} \]

But setting $y = \frac{1}{x}$,
\begin{eqnarray*}
 \int_0^\infty  \log (1+x) \log^{2h-1} x \frac{ \dd x}{x^2 - 1} & = & \int_0^1
\log (1+x) \log^{2h-1} x \frac{ \dd x}{x^2 - 1} \\
&&+ \int_0^1  \log \left (
\frac{1+y}{y}\right) \log^{2h-1} y \frac{ \dd y}{y^2 - 1} \\
& = & 2 \int_0^1
\log (1+y) \log^{2h-1} y\frac{ \dd y}{y^2 - 1}  + \int_0^1 \log^{2h} y \frac{
\dd y}{1 - y^2 } . \end{eqnarray*}

Then the Mahler measure is
\[= \sum_{h=1}^n  \frac{s_{n-h}(2^2, \dots, (2n-2)^2)}{(2n-1)!} 2^{2h}
\pi^{2n+1-2h} \]
\[\left ( \int_0^1  \log (1+y) \log^{2h-1} y\frac{ \dd y}{y^2 - 1}
+ \int_0^1 \log^{2h} y \frac{\dd y}{1-y^2} \right ).\] If we apply  Lemma
\ref{lemz} and Proposition \ref{prop11}, we obtain \[= \sum_{h=1}^n
\frac{s_{n-h}(2^2, \dots, (2n-2)^2)}{(2n-1)!} 2^{2h} \pi^{2n+1-2h} \left (
(2h-1)! (2h+1) \frac{ 2^{2h+1} - 1}{2^{2h+2}} \zeta(2h+1) \right.\] \[+
(-1)^{h-1} \frac{\log 2 (2^{2h}-1) B_{2h} \pi^{2h}}{4 h} \] \[\left. +
\frac{(2h-1)!}{2} \sum_{k=1}^{h-1}
\frac{(2^{2h+1-2k}-1)(2^{2k-1}-1)}{2^{2h-2k}} \frac{(-1)^{k-1} B_{2k} }{
(2k)!}\pi^{2k} \zeta(2h+1-2k)  \right ). \] Finally, by applying equality
(\ref{eqB}) from the Appendix and changing the order of the sums (and setting $t
= h-k$): \[= \frac{\pi^{2n+1}}{2} \log 2 + \sum_{h=1}^n  \frac{s_{n-h}(2^2,
\dots, (2n-2)^2)}{(2n-1)!}   (2h)! \frac{ 2^{2h+1} - 1}{4}
\pi^{2n+1-2h}\zeta(2h+1)\] \[+ \sum_{t=1}^n \frac{(2t)!(2^{2t+1}-1)}{4(2n-1)!}
\] \[\left(\sum_{k=0}^{n-t}  s_{n-t-k}(2^2, \dots, (2n-2)^2) \binom{2(k+t)}{2l}
(-1)^{k-1}\frac{2^{2k}(2^{2k-1}-1)}{k+t}  B_{2k} \right)\pi^{2n+1-2t}
\zeta(2t+1).\]

For the odd  case we need the formula when the parameter $\alpha$ is purely imaginary,
\[\pi m(1+\ii \alpha x +(1-\ii \alpha) y  ) = \frac{\pi}{4} \log\left | \alpha^2 +1 \right | + \im \left ( \Li_2\left ( \ii \left | \alpha \right | \right ) \right )    \]
where $\alpha \in \Rset$.

\[ \pi^{2n+2} m\left( 1+  \left( \frac{1-x_1}{1+x_1}\right) \dots \left( \frac{1-x_{2n+1}}{1+x_{2n+1}}\right) x + \left(1 -  \left( \frac{1-x_1}{1+x_1}\right) \dots \left( \frac{1-x_{2n+1}}{1+x_{2n+1}}\right) \right) y
 \right)\]
\[= \pi 2^{2n+1} \sum_{h=0}^{n} b_{n,h} \left(\frac{\pi}{2}\right)^{2n-2h} \frac{1}{2}\int_{-\infty}^\infty  m(P_{(-1)^n \ii x}) \log^{2h} |x| \frac{ \dd x}{x^2 + 1} \]
\[= 2^{2n+1} \sum_{h=0}^{n} b_{n,h} \left(\frac{\pi}{2}\right)^{2n-2h} \int_{0}^\infty \left (  \frac{\pi}{4} \log (1+x^2)  + \im \left ( \Li_2\left ( \ii x \right ) \right ) \right) \log^{2h} x  \frac{\dd x}{x^2+1} \]

We will now apply Propositions \ref{log2} and \ref{li2},
\[ = \sum_{h=0}^{n} \frac{s_{n-h}(1^2, \dots, (2n-1)^2)}{(2n)!} 2^{2h+1}\pi^{2n-2h} \left( (-1)^h E_{2h} \left( \frac{\pi}{2}\right)^{2h+2}\log 2\right.\]
\[ \left .  +\sum_{l=1 }^{h} \frac{(2h)!} { (2h-2l)!}  (-1)^{h-l} E_{2h-2l} \left( \frac{\pi}{2}\right)^{2h-2l+2} \left( \frac{2^{2l+1}- 1}{2^{2l+1}}\right) \zeta(2l+1)\right.\]
\[ \left. + \sum_{l=0}^{h} B_{2l} \frac{(2h)!}{(2l)!} (2^{2l-1}-1) (-1)^{l+1}
\pi^{2l} (h-l+1) \left( \frac{2^{2h+3-2l} - 1}{2^{2h+1}}
\right)\zeta(2h+3-2l)\right).\]

Applying equation (\ref{eqE}) from the Appendix
\[ = \frac{\pi^{2n+2}}{2} \log 2  + \sum_{h=0}^{n} \frac{s_{n-h}(1^2, \dots, (2n-1)^2)}{(2n)!} 2^{2h+1}\pi^{2n-2h} \]
\[ \left ( \sum_{l=1}^{h} \frac{(2h)!} { (2h-2l)!}  (-1)^{h-l} E_{2h-2l} \left( \frac{\pi}{2}\right)^{2h-2l+2} \left( \frac{2^{2l+1}- 1}{2^{2l+1}}\right) \zeta(2l+1)\right.\]
\[ \left.+  \sum_{l=0}^{h} B_{2l} \frac{(2h)!}{(2l)!} (2^{2l-1}-1) (-1)^{l+1}
\pi^{2l} (h-l+1) \left( \frac{2^{2h+3-2l} - 1}{2^{2h+1}} \right)\zeta(2h+3-2l)
\right).\]

Let us observe the following term carefully,
\[\sum_{h=0}^{n} \frac{s_{n-h}(1^2, \dots, (2n-1)^2)}{(2n)!} 2^{2h+1}\pi^{2n-2h} \]
\[\sum_{l=0}^{h} B_{2l} \frac{(2h)!}{(2l)!} (2^{2l-1}-1) (-1)^{l+1} \pi^{2l}
(h-l+1) \left( \frac{2^{2h+3-2l} - 1}{2^{2h+1}} \right)\zeta(2h+3-2l).\]

Set $s=h-l$,
\[=\frac{1}{(2n)!}\sum_{s=0}^{n} \pi^{2n-2s} \zeta(2s+3) (s+1) (2s)! ( 2^{2s+3} - 1) \]
\[\sum_{l=0}^{n-s} s_{n-s-l}(1^2, \dots, (2n-1)^2) B_{2l} \binom{2(s+l)}{2l}
(2^{2l-1}-1) (-1)^{l+1}.\] By Theorem \ref{thm18} from the Appendix,
\[=\frac{1}{(2n)!}\sum_{s=0}^{n} \pi^{2n-2s} \zeta(2s+3) (s+1) (2s)! ( 2^{2s+3} - 1) \frac{2s+1}{2(2n+1)} s_{n-s}(2^2,\dots,(2n)^2)\]
\[=\frac{1}{(2n+1)!}\sum_{s=0}^{n} \pi^{2n-2s} \zeta(2s+3) (2s+2)! \frac{
2^{2s+3} - 1}{4}  s_{n-s}(2^2,\dots,(2n)^2).\]

Finally, the Mahler measure is
\[ = \frac{\pi^{2n+2}}{2} \log 2  + \frac{1}{(2n+1)!}\sum_{s=0}^{n}   \frac{(2s+2)!( 2^{2s+3} - 1)}{4}  s_{n-s}(2^2,\dots,(2n)^2)\pi^{2n-2s} \zeta(2s+3)\]
\[+ \sum_{l=1 }^{n} \frac{(2l)! (2^{2l+1}-1)}{4 (2n)!} \left(
\sum_{h=0}^{n-l}s_{n-l-h}(1^2, \dots, (2n-1)^2) \binom{2(h+l)}{2l} (-1)^{h}
E_{2h}\right) \pi^{2n-2l+2} \zeta(2l+1).\] Let us also add , that with the help
of Proposition \ref{thmBE} the above equation may be written in terms of
Bernoulli numbers instead of Euler numbers.

\section{ Appendix-- Some identities involving Bernoulli and Euler numbers}

The main result of this section is a collection of identities involving Bernoulli numbers and symmetric functions, which can be deduced from the explicit form of the polynomials $P_k$ and their behavior as it was studied in Section 5. In add ition to those identities, and for completeness, we also mention some other properties of Bernoulli and Euler numbers that have been used in order to simplify the final form of the equations of Theorem 1.

We begin by explicitly computing the polynomials $P_k$:

\begin{prop} \label{propP} We have the following:
\begin{equation}
P_k(x) =- \frac{2}{k+1} \sum_{h=0}^k B_h \binom{k+1}{h} (2^{h-1}-1) \ii ^ h
x^{k+1-h}. \end{equation}
\end{prop}

\pf It is clear that the equation is true for $k=0, 1$.  We will prove that the properties of Lemma \ref{easyprop} hold. But these properties are straightforward except for {\it 4}. Then it is enough to verify property {\it 4}.
\[ -\frac{2l+1}{2 \, \ii^{2l+1}} P_{2l}(\ii) = \sum_{h=0}^{2l} B_h \binom{2l+1}{h} (2^{h-1}-1) \]
Thus, it suffices to prove that
\[ 0 \stackrel{?}{=} \sum_{h=0}^{2l} B_h \binom{2l+1}{h}(2^{h-1}-1)\qquad
\mbox{for}\,\,l>0.\] Using the well known identity:
\begin{equation}
\sum_{s=0}^k \binom{k+1}{s} B_s = 0
\end{equation}
for $k=2l$, we conclude that we only need to prove,
\[ 0 \stackrel{?}{=} \sum_{h=0}^{2l+1} B_h \binom{2l+1}{h}2^{h-1} \qquad \mathrm{for}\,\, l>0\]
since
\[B_{2j+1} = 0 \quad j=1,2,\dots \]
but that is true, because of this other well known identity
\begin{equation}\label{eq:39}
(1-2^{k-1})B_k = \sum_{s=0}^k 2^{s-1} \binom{k}{s} B_s \quad \mbox{for} \, \,
n>1. \end{equation} \qed

Let us mention the following technical consequence that will be used later.
\begin{cor} We have the following special values:
\begin{equation}
P_{2l-1}(\ii) = (-1)^{l}\frac{2^{2l}-1}{l}B_{2l}.
\end{equation}
\end{cor}
\pf
\[P_{2l-1}(\ii) = -\frac{1}{l} \sum_{h=0}^{2l-1} B_h \binom{2l}{h} (2^{h-1}-1) (-1)^l =
\frac{(-1)^{l+1}}{l} \sum_{h=0}^{2l-1} B_h \binom{2l}{h} 2^{h-1}\]
\[=\frac{(-1)^{l+1}}{l} ((1-2^{2l-1})B_{2l}-2^{2l-1} B_{2l}) =\frac{(-1)^{l+1}}{l} (1-2^{2l})B_{2l}\]
because of equation (\ref{eq:39}). \qed

In fact,
\begin{cor} We have:
\begin{equation}
P_{k}(x) = \frac{2 \ii^{k+1}}{k+1} \left(B_{k+1}\left(\frac{x}{\ii}\right) -
2^{k} B_{k+1}\left(\frac{x}{2\ii}\right)
\right)+\frac{(2^{k+1}-2)\ii^{k+1}}{k+1} B_{k+1} \end{equation} where $B_k(x)$
is the Bernoulli polynomial. \end{cor}

We are now ready to prove the main Theorem of this section.
\begin{thm} \label{thm18} We have the following identities:

For $1\leq l \leq n$:
\[s_{n-l}(1^2, \dots, (2n-1)^2) \]
\[= n \sum_{s=0}^{n-l} s_{n-l-s}(2^2, \dots, (2n-2)^2)  \frac{1}{l+s} B_{2s}
\binom{2(l+s)}{2s} (2^{2s}-2) (-1)^{s+1}.\]

For $1 \leq n$:
\[\left(\frac{(2n)!}{2^n n!} \right)^2 = 2n \sum_{s=1}^{n} s_{n-s}(2^2, \dots,
(2n-2)^2)  \frac{1}{s} B_{2s} (2^{2s}-1) (-1)^{s+1}. \]

For $0 \leq l \leq n$:
\[(2l+1) s_{n-l}(2^2, \dots, (2n)^2) \]
\[=(2n+1)\sum_{s=0}^{n-l} s_{n-l-s}(1^2, \dots, (2n-1)^2) B_{2s}
\binom{2(l+s)}{2s} (2^{2s}-2) (-1)^{s+1}.\]
\end{thm}

\pf By Lemma \ref{lemma13} and Theorem \ref{thm14} we have
\[\frac{(x^2+1^2) \dots (x^2+(2n-1)^2)}{(2n)!} \]
\[=   \sum_{h=1}^{n} \frac{s_{n-h}(2^2, \dots, (2n-2)^2)}{(2n-1)!}  \left( -
\frac{2}{2h} \sum_{j=0}^{2h-1} B_j \binom{2h}{j} (2^{j-1}-1) \ii ^ j x^{2h-j} -
(-1)^{h}\frac{2^{2h}-1}{h}B_{2h}\right).\] Set $j=2s$, then the first term in
the difference is \[ \sum_{h=1}^{n} \frac{s_{n-h}(2^2, \dots,
(2n-2)^2)}{(2n-1)!} \left( - \frac{1}{h} \sum_{s=0}^{h-1} B_{2s} \binom{2h}{2s}
(2^{2s-1}-1) (-1)^s x^{2h-2s} \right).\] Now set $l=h-s$, \[=   \sum_{h=1}^{n}
\frac{s_{n-h}(2^2, \dots, (2n-2)^2)}{(2n-1)!}  \left( \frac{1}{h} \sum_{l=1}^{h}
B_{2(h-l)} \binom{2h}{2l} (2^{2(h-l)-1}-1) (-1)^{h-l +1} x^{2l} \right)\] \[=
\sum_{l=1}^n \left(\sum_{h=l}^{n} \frac{s_{n-h}(2^2, \dots, (2n-2)^2)}{(2n-1)!}
\frac{1}{h} B_{2(h-l)} \binom{2h}{2l} (2^{2(h-l)-1}-1) (-1)^{h-l+1}\right)
x^{2l}.\]

Comparing coefficients we get
\[\frac{s_{n-l}(1^2, \dots, (2n-1)^2)}{(2n)!} \]
\[=\sum_{h=l}^{n}
\frac{s_{n-h}(2^2, \dots, (2n-2)^2)}{(2n-1)!}  \frac{1}{h} B_{2(h-l)}
\binom{2h}{2l} (2^{2(h-l)-1}-1) (-1)^{h-l+1}.\]
Thus\[s_{n-l}(1^2, \dots,
(2n-1)^2)\]
\[= n \sum_{s=0}^{n-l} s_{n-l-s}(2^2, \dots, (2n-2)^2)  \frac{1}{l+s} B_{2s}
\binom{2(l+s)}{2s} (2^{2s}-2) (-1)^{s+1}.\]

The second equality is obtained by comparing the independent coefficients.

For the third equality, we do a similar process:
\[\frac{x(x^2+2^2)\dots(x^2+(2n)^2)}{(2n+1)!} \]
\[=   \sum_{h=0}^{n} \frac{s_{n-h}(1^2, \dots, (2n-1)^2)}{(2n)!} \left(-
\frac{2}{2h+1} \sum_{j=0}^{2h} B_j \binom{2h+1}{j} (2^{j-1}-1) \ii ^ j
x^{2h+1-j}\right)\]
\[ =\sum_{h=0}^{n} \frac{s_{n-h}(1^2, \dots,
(2n-1)^2)}{(2n)!} \left(- \frac{2}{2h+1} \sum_{s=0}^{h} B_{2s} \binom{2h+1}{2s}
(2^{2s-1}-1) (-1)^s x^{2h+1-2s}\right).\] Now set $l=h-s$, \[ =\sum_{h=0}^{n}
\frac{s_{n-h}(1^2, \dots, (2n-1)^2)}{(2n)!} \left( \frac{2}{2h+1} \sum_{l=0}^{h}
B_{2(h-l)} \binom{2h+1}{2l+1} (2^{2(h-l)-1}-1) (-1)^{h-l+1} x^{2l+1}\right)\]
\[=\sum_{l=0}^n \left(\sum_{h=l}^{n} \frac{s_{n-h}(1^2, \dots,
(2n-1)^2)}{(2n)!}\frac{2}{2h+1}B_{2(h-l)} \binom{2h+1}{2l+1} (2^{2(h-l)-1}-1)
(-1)^{h-l+1} \right)x^{2l+1}.\]

Comparing coefficients we get
\[ \frac{s_{n-l}(2^2, \dots, (2n)^2)}{(2n+1)!} \]
\[=\sum_{h=l}^{n}
\frac{s_{n-h}(1^2, \dots, (2n-1)^2)}{(2n)!}\frac{2}{2h+1}B_{2(h-l)}
\binom{2h+1}{2l+1} (2^{2(h-l)-1}-1) (-1)^{h-l+1}.\]
Thus
\[ (2l+1) s_{n-l}(2^2, \dots,
(2n)^2) \]
\[=(2n+1)\sum_{s=0}^{n-l} s_{n-l-s}(1^2, \dots, (2n-1)^2) B_{2s}
\binom{2(l+s)}{2s} (2^{2s}-2) (-1)^{s+1}.\] \qed

The next result illuminates the last formula of Theorem 1.
\begin{prop} \label{thmBE} We have
\[n \sum_{s=0}^{n-l} s_{n-l-s}(2^2, \dots, (2n-2)^2)  \frac{1}{l+s} B_{2s} \binom{2(l+s)}{2s} 2^{2s}(2^{2s}-2) (-1)^{s+1}\]
\begin{equation} \label{crazy}
=\sum_{k=l}^{n}  (-1)^{k+l}\binom{2k}{2l} s_{n-k}(1^2, \dots, (2n-1)^2)
E_{2(k-l)}. \end{equation}
\end{prop}

\pf By Proposition \ref{propsym},
\[ n \sum_{s=0}^{n-l} s_{n-l-s}(2^2, \dots, (2n-2)^2)  \frac{1}{l+s} B_{2s} \binom{2(l+s)}{2s} 2^{2s}(2^{2s}-2) (-1)^{s+1}\]
\[= \sum_{s=0}^{n-l} \frac{(-1)^{l+s}}{2} \]
\[\sum_{k=l+s}^n (-1)^k
\binom{2k}{2s+2l-1} s_{n-k}(1^2, \dots, (2n-1)^2)  \frac{1}{l+s} B_{2s}
\binom{2(l+s)}{2s} 2^{2s}(2^{2s}-2) (-1)^{s+1}\] \[= \sum_{s=0}^{n-l}
\sum_{k=l+s}^n \frac{(-1)^{k+l+1} }{2k+1}\binom{2k+1}{2l} s_{n-k}(1^2, \dots,
(2n-1)^2)   B_{2s} \binom{2(k-l)+1}{2s} 2^{2s}(2^{2s}-2).\] Changing the order
of the sums, \[= \sum_{k=l}^{n}   \frac{(-1)^{k+l+1} }{2k+1}\binom{2k+1}{2l}
s_{n-k}(1^2, \dots, (2n-1)^2)  \sum_{s=0}^{k-l} B_{2s} \binom{2(k-l)+1}{2s}
2^{2s}(2^{2s}-2).\]

Now observe that
 \[\sum_{s=0}^{k-l} B_{2s} \binom{2(k-l)+1}{2s} 2^{2s}(2^{2s}-2) =
\sum_{m=0}^{2(k-l)+1}  B_m \binom{2(k-l)+1}{m} 2^m(2^m-2).\] By equation
(\ref{eq:39}), \[ = \sum_{m=0}^{2(k-l)+1}  B_m \binom{2(k-l)+1}{m} 2^{2m}=
4^{2(k-l)+1} B_{2(k-l)+1}\left(\frac{1}{4}\right)\] for $k-l>0$ and $=-1$
otherwise.

Now use the following identity,
\[ 2^{2n}B_n\left(\frac{1}{4} \right)  =(2- 2^n) B_n - n E_{n-1},\]
which can be found, for instance, as equation (23.1.27) in \cite{AS}.



Then, we get:
\[\sum_{s=0}^{k-l} B_{2s} \binom{2(k-l)+1}{2s} 2^{2s}(2^{2s}-2) =  - (2(k-l)+1)
E_{2(k-l)}.\]

Therefore,
\[ -\sum_{k=l}^{n}   \frac{(-1)^{k+l+1} }{2k+1}\binom{2k+1}{2l}
s_{n-k}(1^2, \dots, (2n-1)^2) (2(k-l)+1) E_{2(k-l)}\]
\[ = \sum_{k=l}^{n}  (-1)^{k+l}\binom{2k}{2l}
s_{n-k}(1^2, \dots, (2n-1)^2)  E_{2(k-l)}.\] \qed

We would like to finish by stating a few basic equalities that can be proved by induction:
\begin{prop} \label{propBEsimple} We have:
\begin{equation} \label{eqE}
 \sum_{h=0}^n s_{n-h}(1^2, \dots, (2n-1)^2) (-1)^h E_{2h} = (2n)!
 \end{equation}
\begin{equation}
 \sum_{h=0}^n s_{n-h}(1^2, \dots, (2n-1)^2) (-1)^{h+1} E_{2(h+1)} = (2n + 1)!
 \end{equation}
 \begin{equation} \label{eqB}
 \sum_{h=1}^n s_{n-h}(2^2, \dots, (2n-2)^2) (-1)^{h+1}\frac{2^{2h}(2^{2h}-1)}{h} B_{2h} = 2 (2n - 1)!
 \end{equation}
\end{prop}

\section{Concluding remarks}
In conclusion, the Mahler measure of these three families of $n$-variable polynomials can be computed explicitly as some linear combination of special values of zeta functions, the L-series on the Dirichlet character of conductor 4, (and $\mathcal{L}_{3,2h+1}(\ii,\ii)$ for the second family). It remains to relate $\mathcal{L}_{3,2h+1}(\ii,\ii)$ to L-series and perhaps zeta functions, which would simplify formula (\ref{eq4}).

In some cases the coefficients of these formulas are related to Bernoulli
numbers. It should be remarked that the results of Theorem \ref{thm18} and
Proposition \ref{propBEsimple} suggest that there should be a simpler expression
for formulas of the kind of Theorem \ref{thmBE}, and that might allow to find
better expressions for the formulas of case (iii) (equations (\ref{eq5}) and
(\ref{eq6})), for instance.

Finally and most importantly, it would be interesting to find different families, perhaps, by add ing new variables by using other forms of fractional transformations or other rational functions.

\bigskip
\noindent{\bf Acknowledgements}

It is my pleasure to thank Fernando Rodriguez-Villegas for several helpful
discussions and his constant encouragement. I am also grateful to David Boyd
for his useful comments, and to the Referee, whose suggestions have
improved the exposition of this article.

This work was developed while I was receiving financial support from the Harrington fellowship and John Tate to whom I am deeply grateful.

 \end{document}